\pdfoutput=1
\RequirePackage{ifpdf}
\ifpdf 
\documentclass[pdftex]{sigma}
\else
\documentclass{sigma}
\fi

\numberwithin{equation}{section}

\newtheorem{Theorem}{Theorem}[section]
\newtheorem*{Theorem*}{Theorem}
\newtheorem{Corollary}[Theorem]{Corollary}
\newtheorem{Lemma}[Theorem]{Lemma}
\newtheorem{Proposition}[Theorem]{Proposition}

\theoremstyle{definition}
\newtheorem{Definition}[Theorem]{Definition}

\newtheorem{Example}[Theorem]{Example}
\newtheorem{Remark}[Theorem]{Remark}

\usepackage[all]{xy}
\usepackage{tikz-cd}

\newcommand{\cF}{\mathcal{F}}
\newcommand{\cG}{\mathcal{G}}
\newcommand{\ZZ}{\ensuremath{\mathbb Z}}
\newcommand{\RR}{\ensuremath{\mathbb R}}
\newcommand{\g}{\ensuremath{\mathfrak{g}}}
\newcommand{\pd}[1]{\frac{\partial}{\partial #1}} 
\newcommand{\vX}{\mathfrak{X}}
\newcommand{\cD}{\mathcal{D}}
\newcommand{\bs}{\mathbf{s}}
\newcommand{\bt}{\mathbf{t}}
\newcommand{\w}[1]{\widetilde{#1}}

\begin{document}

\allowdisplaybreaks

\newcommand{\arXivNumber}{2506.19509}

\renewcommand{\thefootnote}{}

\renewcommand{\PaperNumber}{014}

\FirstPageHeading

\ShortArticleName{The Linearizability of Singular Foliations Is a Morita Invariant}

\ArticleName{The Linearizability of Singular Foliations\\ Is a Morita Invariant\footnote{This paper is a~contribution to the Special Issue on Interactions of Poisson Geometry, Lie Theory and Symmetry in honor of Rui Loja Fernandes for his 60th birthday. The~full collection is available at \href{https://sigma-journal.com/Fernandes.html}{https://sigma-journal.com/Fernandes.html}}}

\Author{Marco ZAMBON}

\AuthorNameForHeading{M.~Zambon}

\Address{KU Leuven, Department of Mathematics, Celestijnenlaan 200B box 2400, 3001 Leuven, Belgium}
\Email{\mail{marco.zambon@kuleuven.be}}
\URLaddress{\url{https://marcozambonmath.github.io/}}

\ArticleDates{Received September 25, 2025, in final form February 05, 2026; Published online February 17, 2026}

\Abstract{Hausdorff Morita equivalence is an equivalence relation on singular foliations, which induces a bijection between their leaves. Our main statement is that linearizability along a leaf is invariant under Hausdorff Morita equivalence. The proof relies on a characterization of tubular neighborhood embeddings using Euler-like vector fields.}


\Classification{53C12; 53D17}

\renewcommand{\thefootnote}{\arabic{footnote}}
\setcounter{footnote}{0}

\section{Introduction}

This note is concerned with \emph{singular foliations}, understood as a module of vector fields as in \cite{AndrSk}, following the work of Stefan and Sussman.
{A singular foliation $\cF$ on a manifold $M$} gives rise to a decomposition of $M$ into immersed submanifolds tangent to $\cF$, called leaves.
For instance, any action of a Lie group $G$ on a manifold $M$ gives rise to such a singular foliation, namely the $C^{\infty}_c(M)$-span of the vector fields generating the action. The underlying decomposition of $M$ is by the orbits of the $G$-action, if $G$ is connected.

Whenever $L$ is an embedded leaf of a singular foliation $\cF$, by linearizing the vector fields in~$\cF$, one obtains a singular foliation $\cF_{\rm lin}$ on the normal bundle $\nu L$, called linearized foliation. The \emph{linearization} question asks when the singular foliations $\cF$ and $\cF_{\rm lin}$ are isomorphic, in a~neighborhood of~$L$. When $L$ is a point, sufficient criteria for linearization
were given by~Cerveau~\cite{CerveauSing}.\looseness=1

\emph{Hausdorff Morita equivalence} \cite{HME} is a notion of equivalence for singular foliations, weaker than isomorphism, which among other things preserves the geometry transverse to the leaves, and induces a homeomorphism of leaf spaces. Our main result is that it preserves linearizability.\looseness=1

 \begin{theorem*}[Theorem~\ref{thm:main}]
Let $(M_1,\cF_1)$ and $(M_2,\cF_2)$ be Hausdorff Morita equivalent singular foliations,
and $L_1\subset M_1$, $L_2\subset M_2$ corresponding embedded leaves.
Then $\cF_1$ is linearizable around $L_1$ if and only if $\cF_2$ is linearizable around $L_2$.
\end{theorem*}

We remark that the dimensions of $L_1$ and $L_2$ are generally different,
hence this theorem allows to reduce the linearization question to leaves of smaller dimension, which in certain cases are simply points.

By the very definition of Hausdorff Morita equivalence, the theorem follows immediately from the following proposition.

\begin{proposition*}[Proposition~\ref{prop:pss}]
 Let $\pi\colon P\to M$ be a surjective submersions with connected fibers. Let $\cF$ be a singular foliation on $M$ and $L\subset M$ an embedded leaf.
 Then
 $\cF$ is linearizable around $L$ if and only if the pullback foliation $\pi^{-1}\cF$ is linearizable around $\pi^{-1}(L)$.
\end{proposition*}

 The main tool we use to prove this result is the correspondence between \emph{tubular neighborhood embeddings}
 and \emph{Euler-like vector fields} \cite{DefSpacesEulerLike,SplittingThmEulerlike,HigsonReza}. The precise relation between linearization and Euler-like vector fields is given in Proposition~\ref{prop:Cerveau}, which extends a result of \cite[Section~8]{CerveauSing} for points.
With this tool at hand, the direct implication in the proposition is easy to prove, by lifting Euler-like vector fields. The converse implication requires more work, since a given Euler-like vector field is not necessarily $\pi$-projectable.

In the last section of the paper,
 we present examples. Further, using the above theorem, we~obtain a condition under which the linearizability of a singular foliation $\cF$ along a leaf is implied by the linearizability of a simpler foliation, namely the restriction of $\cF$ to a slice.

\begin{corollary*}[Corollary~\ref{cor:slicelin}]
Let $\cF$ be a singular foliation arising from a Hausdorff, source connected, Lie groupoid $\cG$.
Let $L$ be an embedded leaf, $N$ a slice. Suppose the restricted Lie groupoid $\cG_N:=\bt^{-1}(N)\cap \bs^{-1}(N)$ is source connected.
Then, whenever $\bigl(N,\iota_N^{-1}\cF\bigr)$ is linearizable {around the point $N\cap L$}, the singular foliation $\cF$ is linearizable {around $L$}.
 \end{corollary*}

 \section{Singular foliations and their linearization}

We present some facts about singular foliations, their linearization, and their Hausdorff Morita equivalence. We also recall a theorem about tubular neighborhood embeddings.

 \subsection{Singular foliations}\label{sec:singfol}
 \label{subsec:sf}

We recall the notion of singular foliation from \cite{AndrSk}, see also the comprehensive monograph \cite{singfolnotes}.

 \begin{Definition}
A \textit{singular foliation} on a manifold $M$ is a $C^{\infty}(M)$-submodule $\cF$
{of the compactly supported vector fields}
$\vX_c(M)$, involutive with respect to the Lie bracket and locally finitely generated.
A \textit{foliated manifold} is a manifold with a singular foliation.
\end{Definition}

\begin{Remark}\label{rem:fhat}
 For any open set $U\subset M$,
consider the following modules:
\begin{align*}
&\cF|_U:=\iota_U^{-1}\cF:=\{ X|_U: X\in \cF \text{ and } \text{supp}(X)\subset U\},\\
&\widehat{\iota_U^{-1}\cF}:=\bigl\{ X\in \vX(U) : fX\in \iota_U^{-1}\cF \ \text{ for all }f\in C^{\infty}_c(U)\bigr\}.
\end{align*}
Notice that if $Z\in \cF$, its restriction $Z|_U$ is not compactly supported in general, and \smash{$Z|_U\in \widehat{\iota^{-1}_{U}\cF}$}.
We say that $\cF$ is \emph{locally finitely generated} if for every point of $M$ there is an open neighborhood~$U$ and finitely many \smash{$X_1,\dots,X_n\in \raisebox{-1pt}{$\widehat{\iota_U^{-1}\cF}$}$}
{such that \smash{$\iota^{-1}_{U}\cF={\rm Span}_{C^{\infty}_c(U)}\{Y_i\}$}.}
 In that case, \smash{$\raisebox{-1pt}{$\widehat{\iota^{-1}_{U}\cF}$}={\rm Span}_{C^{\infty}(U)}\{Y_i\}$}, see \cite[Example~3\,(iii)]{AZ6}. When $U=M$, we obtain the global hull of $\cF$:
\[
\widehat{\cF}:=\{ X\in \vX(M) : fX\in \cF \ \text{ for all }f\in C^{\infty}_c(M)\}.
\]
\end{Remark}

 \begin{Remark}
 A singular foliation on a manifold $M$ can be equivalently regarded as an involutive, locally finitely generated subsheaf of
the sheaf of $C^{\infty}$-modules on $M$ given by the smooth vector fields. We will not use this point of view in this paper.
\end{Remark}

 By the Stefan--Sussmann theorem, a singular foliation induces a partition of the manifold into immersed submanifolds, called \emph{leaves}. For results on the structure of singular foliations nearby a given leaf, see
\cite{CamilleSR,francis2024singularfoliationstangentgiven,NbhdSingLeaf}.

\begin{Example} Any Lie algebroid $A$ induces a singular foliation on its base manifold $M$, namely $\rho(\Gamma_c(A))$, where $\rho\colon A\to TM$ denotes the anchor map.
In particular, any Lie groupoid over $M$ induces a singular foliation (via its Lie algebroid).
 \end{Example}

 \begin{Definition}\label{rem:pulblackfol}
 Let $(M,\cF)$ be a foliated manifold and $\pi\colon P\to M$ a submersion.
The {\it pullback foliation}~$\pi^{-1}\cF$
\cite[Proposition 1.10]{AndrSk}
is the $C^{\infty}_c(P)$-span of the set of projectable vector fields on~$P$ which project to elements of~$\cF$.
\end{Definition}

In this definition, one needs to take the $C^{\infty}_c(P)$-span (rather than the $C^{\infty}(P)$-span) because projectable vector fields are usually not compactly supported.
In \cite[Proposition~1.10]{AndrSk},
 the pullback foliation is defined for any map transverse to $\cF$. We will need the case of the inclusion $\iota_S \colon S\to M$ of a submanifold transverse to the foliation: in that case,
 \[
 \iota_S^{-1}\cF={\rm Span}_{C^{\infty}_c(S)}\{X|_S: X\in \cF \text{ is tangent to $S$}\}.
 \]

We now consider flows. In \cite[Proposition~1.6]{AndrSk} (see also \cite[Proposition~2.3]{OriAlfonso}), it is shown that for any singular foliation $\cF$ and $X\in \cF$, the time-1 flow of $X$ preserves $\cF$. We need an~extension of that result to infinitesimal symmetries, i.e., vector fields $X$ satisfying ${[X,\cF]\subset \cF}$. For compactly supported $X$ this appears in \cite[Proposition~2.1.3]{AlfonsoThesis}, and in the form presented here it appears in \cite[Proposition~1.7.10]{singfolnotes} (see also \cite[Proposition~1.3]{symsingfol}).

\begin{Proposition}\label{prop:OriAlfcompl}
Let $(M,\cF)$ be a foliated manifold and $X$ a vector field on $M$ such that $[X,\cF]\subset \cF$. Then the time-$t$ flow $\phi_t$ of $X$, where it exists, satisfies $(\phi_t)_*\cF\subset \cF$.
\end{Proposition}

We will need a technical lemma.
\begin{Lemma}\label{lem:average}
Let $(M,\cF)$ be a foliated manifold, $\{X_t\}_{t\in [0,1]}$ a smooth family of vector fields lying in $\cF$. Assume that there exists a cover $\{U_{\alpha}\}$ of $M$ by open subsets so that on each $U_{\alpha}$ the singular foliation $\cF$ admits a finite number of generators $Y_1^{\alpha},\dots,Y^{\alpha}_{k_{\alpha}}\in \vX(U_{\alpha})$
with the following property:
\[
X_t|_{U_{\alpha}}=\sum_i f_i^{t,\alpha}Y_i^{\alpha},
\]
where the coefficients $f_i^{t,\alpha}\in C^{\infty}(U_{\alpha})$ depend smoothly on $t$. Then the vector field
\[
 X:=\int_0^1 X_t \, {\rm d}t
\]
lies in the global hull $\widehat{\cF}$.
\end{Lemma}
\begin{proof}
Take a cover $\{U_{\alpha}\}$, generators $Y_i^{\alpha}$ and functions \smash{$f_i^{t,\alpha}$} as above.
Since
\[
\biggl(\int_0^1 f_i^{t,\alpha} {\rm d}t\biggr)\in C^{\infty}(U_{\alpha}),
\]
we can write
\[
X|_{U_{\alpha}}=\sum_i \biggl(\int_0^1 f_i^{t,\alpha} {\rm d}t\biggr)Y_i^{\alpha},
\]
which lies in \smash{$\widehat{\iota^{-1}_{U_{\alpha}}\cF}$} by Remark~\ref{rem:fhat}.

 Fix a partition of unity $\{\rho_{\alpha}\}$ subordinate to the above cover of $M$.
Take $f\in C^{\infty}_c(M)$. Notice that $\{\alpha: {\rm supp}(\rho_{\alpha})\cap {\rm supp}(f)\neq \varnothing\}$ is a finite set, as a consequence of the partition of unity being locally finite (see \cite[proof of Proposition~1]{AZ6}).
We have
\[
fX=\sum_{\alpha}f\rho_{\alpha}X.
\]
This is a~finite sum, and each summand lies in $\cF$ because we showed that \smash{$X|_{U_{\alpha}}\in \widehat{\iota^{-1}_{U_{\alpha}}\cF}$} and $ f\rho_{\alpha}\in C^{\infty}_c(M)$ has support in $U_{\alpha}$. Hence $fX\in \cF$. As this holds for all $f\in C^{\infty}_c(M)$, we conclude that \smash{$X\in \widehat{\cF}$}.
\end{proof}

\begin{Remark}
The statement of Lemma~\ref{lem:average} no longer holds if we omit the assumption on the smooth dependence of the \smash{$f_i^{t,\alpha}$} on~$t$.
\end{Remark}

We now look at how the flow of an element of $\cF$ acts on an infinitesimal symmetry of $\cF$.

\begin{Proposition}\label{prop:uptoF}
Let $(M,\cF)$ be a foliated manifold, $Y\in \vX(M)$ such that $[Y,\cF]\subset \cF$.
Fix $X\in \cF$, denote its time-$1$ flow by $\phi^X_1$. Then
\[
\bigl(\phi^X_1\bigr)_*Y-Y\in \widehat{\cF}.
\]
\end{Proposition}

\begin{proof}
Notice that
 \[
 \bigl(\phi^X_1\bigr)_*{Y}-{Y}=\int_0^1\frac{{\rm d}}{{\rm d}t}\bigl(\phi^X_t\bigr)_*{Y} \,{\rm d}t;
 \]
 we would like to apply Lemma~\ref{lem:average} to conclude that
 lies in $\widehat{{\cF}}$.

To apply Lemma~\ref{lem:average}, notice that for all $t\in [0,1]$ the integrand satisfies
 \begin{equation}\label{eq:phiyx}
\frac{{\rm d}}{{\rm d}t}\bigl(\phi^X_t\bigr)_*{Y}=\bigl(\phi^X_t\bigr)_*[Y,X]
\end{equation}
by a standard computation (see for instance \cite[proof of Proposition~2.3]{OriAlfonso}).
This lies in ${\cF}$, since $[{Y},X]\in {\cF}$ by assumption and since
$\phi^X_t$ is an automorphism of $(M,{\cF})$ by \cite[Proposition~1.6]{AndrSk}.

We now check that the smooth family of vector fields \eqref{eq:phiyx} satisfies
the assumptions of Lemma~\ref{lem:average}.
Since ${\rm supp}(X)$ is compact, $\cF$ is finitely generated on an open subset $U$ of $M$ containing ${\rm supp}(X)$, as shown by the argument at the beginning of \cite[proof of Proposition~2.3]{OriAlfonso}. We denote a set of generators by $Y_1,\dots,Y_{k}\in \vX(U)$.

Since $[Y,X]\in \cF$ is supported on $U$, we can write $[Y,X]|_U=\sum_i g_iY_i$ for
$g_i\in C^{\infty}_c(U)$. It~is shown in \cite[proof of Proposition~2.3]{OriAlfonso} that \smash{$\bigl(\phi^X_{t}\bigr)_*Y_i=\sum_jf_j^i(t)Y_j$} for some
 \smash{$f_j^i(t)\in C^{\infty}(U)$} depending smoothly on $t$.
Hence on $U$ the right-hand side of~\eqref{eq:phiyx} equals
\[
\sum_i \bigl(g_i\circ \phi^X_{-t}\bigr)\bigl(\phi^X_{t})_*Y_i=
\sum_{i,j} \bigl(g_i\circ \phi^X_{-t}\bigr)f_j^i(t)Y_j.
\]
As the coefficients
\smash{$\bigl(g_i\circ \phi^X_{-t}\bigr)f_j^i(t)$} depend smoothly on $t$, {and as the integrand \eqref{eq:phiyx} vanishes on $M\setminus {\rm Supp}(X)$,} the assumptions of Lemma~\ref{lem:average} are satisfied.
\end{proof}

\subsection{The linearized foliation}\label{subsec:linfol}

We recall the linearization of a singular foliation along an (embedded) leaf, following \cite[Sections~4.1 and~4.3]{AZ2}.
Consider a singular foliation $(M,\cF)$ and an embedded leaf $L$. Denote by $\nu L:=TM|_L/TL$ the normal bundle of $L$.
Recall that there is a canonical identification
\[
I_L/I^2_L \cong \Gamma(\nu^*L)=C^{\infty}_{\rm lin}(\nu L),\qquad [f] \mapsto {\rm d}f|_L,
\]
where $I_L$
 denotes the functions on $M$ that vanish on $L$, \smash{$C^{\infty}_{\rm lin}(\nu L)$} the fiberwise linear functions on the normal bundle, and
 \smash{$[f]:=\bigl(f \bmod I^2_L\bigr)$} for $f\in I_L$.

 Every vector field $Y$ on $M$ tangent to $L$ gives rise to a fiberwise linear vector field $Y_{\rm lin}$ on $\nu L$, which acts as follows on the fiberwise constant functions and on $C^{\infty}_{\rm lin}(\nu L)\cong I_L/I^2_L$:{\samepage
 \begin{align*}
 \begin{split}
&Y_{\rm lin}(g):= (Y|_L)(g) \qquad\text{for all } g\in C^{\infty}(L),\\
&Y_{\rm lin}[f]:=[Y(f)] \qquad\text{for all } f\in I_L.
\end{split}
\end{align*}
Here, given a function on $L$, abusing notation we denote its pullback to $\nu L$ with the same letter.}

We denote by $\chi_{\rm lin}(\nu L)$ the Lie algebra of fiberwise linear vector fields on $\nu L$, i.e., those which preserve the fiberwise constant functions and the fiberwise linear functions.
The above yields a~linear map
\[
{\rm lin} \colon\ \cF \to \chi_{\rm lin}(\nu L),\qquad Y \mapsto Y_{\rm lin}
\]
which preserves brackets.
Notice that $lin$ factors through a map of $C^{\infty}(L)$-modules
\begin{equation}\label{eq:overlinelin}
\overline{\rm lin}\colon \cF/I_L\cF \to \chi_{\rm lin}(\nu L).
\end{equation}
 \begin{Definition}
The {\it linearization of $\cF$ along $L$} is the singular foliation $\cF_{\rm lin}$ on $\nu L$
generated by the image of ${\rm lin}$, i.e., the $C^{\infty}_c(\nu L)$-span of $\{Y_{\rm lin}:Y \in \cF\}$.
\end{Definition}

\begin{Definition}
We say that {\it $\cF$ is linearizable} about $L$ if
 there exist neighborhoods $V\subset \nu L$ of the zero section and $U\subset M$ of $L$, and a tubular neighborhood embedding (see Definition~\ref{def:tneeulerlike})
 $V\to U$
inducing an isomorphism between
 the singular foliations $\cF_{\rm lin}|_V$ and $\cF|_U$.
\end{Definition}

Not all singular foliations are linearizable. For instance, the singular foliation $\cF$ on the real line generated by $x^2\partial_{x}$ is not linearizable along the leaf $L=\{0\}$: the linearized foliation $\cF_{\rm lin}$ is the zero foliation, for which each point is a leaf.

\begin{Remark}\label{rem:flin}
We have $\cF/I_L\cF=\Gamma_c(A_L)$ for a transitive Lie algebroid $A_L$ over $L$ (see \cite[Section~1.3]{AZ1}), whose anchor we will denote $\rho_{A_L}$.
Similarly, $\chi_{\rm lin}(\nu L)$ agrees with the sections of~${\rm CDO}(\nu L)$, the transitive Lie algebroid whose sections are the covariant differential operators on the vector bundle $\nu L$, and whose anchor $\sigma$ is given by the symbol map. Hence the map $\overline{\rm lin}$ in~\eqref{eq:overlinelin} is the map of sections induced by a morphism of transitive Lie algebroids
\[
\lambda\colon\ A_L\to {\rm CDO}(\nu L)
\]
covering ${\rm Id}_L$. Any such morphism has constant rank by \cite[Theorem~6.5.3]{MackenzieGrdAld}. Consequently, $\{Y_{\rm lin}:Y \in \cF\}=\Gamma_c(E)$ for a wide Lie subalgebroid $E\subset {\rm CDO}(\nu L)$.
\end{Remark}

Let $(M,\cF)$ be a foliated manifold and $p\in M$.
Consider \[\g_p:=\{X\in \cF: X(p)=0 \}/I_p\cF,\] where
 $I_p$ denotes the ideal of smooth functions on $M$ vanishing at $p$.
Then $\g_p$ is a Lie algebra, called \textit{isotropy Lie algebra}, with Lie bracket induced by the one of vector fields. It is the isotropy Lie algebra $\ker(\rho_{A_L})_p$ at $p$ of the Lie algebroid $A_L$, where $L$ denotes the leaf of $\cF$ through $p$.

\begin{Remark}
The restriction of the map $\lambda$
to fibers over $p$ is a Lie algebra morphism
\begin{equation}\label{eq:notinj}
\lambda_p \colon\ \g_p \to \vX_{\rm lin}(\nu_p L)={\rm End}(\nu_p L),\qquad Y \bmod I_p\cF \mapsto Y_{\rm lin}|_{\nu_p L}.
\end{equation}
This map is not injective in general, as can be seen from the example where $M=\RR$, $\cF$ is generated by $x^2\partial_{x}$, and $p=0$.
\end{Remark}

 \subsection{Interlude: Tubular neighborhood embeddings}\label{subsec:tne}
Putting aside singular foliations for a moment, here we recall some material on tubular neighborhood embeddings, following \cite[Section~1]{HigsonReza}.

Let $N$ be a (embedded) submanifold of a manifold $M$, denote by $\nu N=TM|_N/TN$ its normal bundle.
\begin{Definition}\label{def:tneeulerlike}\samepage
\quad
\begin{enumerate}\itemsep=0pt
\item[(i)] A {\it tubular neighborhood embedding} is an embedding $\psi\colon V\to M$, {where $V$ is an open neighborhood of the zero section in $ \nu N$,} such that $\psi|_N={\rm Id}_N$ and the induced map on normal bundles is ${\rm Id}_{\nu N}$.

\item[(ii)] A vector field\footnote{We do not require $X$ to be complete, just as in \cite{HigsonReza}. This is in contrast with \cite{SplittingThmEulerlike}, where the vector field is assumed to be complete, and tubular neighborhood embeddings are defined on the whole of $\nu N$.} $X$ on $M$ is {\it Euler-like for $N$} if $X|_N=0$ and $X_{\rm lin}$ is the Euler vector field~$E$ of $\nu N$.
\end{enumerate}
\end{Definition}

\begin{Remark} \label{rem:charEuler} The Euler vector field $E$ of $\nu N$
 is fiberwise linear, and is the unique vector field satisfying the following properties:
$E(g)=0 \text{ for all } g\in C^{\infty}(N)$,
$E[f] =[f] \text{ for all } f\in I_N$. {From this, one sees that the Euler vector field commutes with all fiberwise linear vector fields.}
\end{Remark}

Clearly if $\psi\colon V\to M$ is a tubular neighborhood embedding, then pushing forward the Euler vector field of $\nu N$ one obtains a Euler-like vector field\footnote{{This Euler-like vector field is not complete, unless $V=\nu N$.}} for $N$, defined on $\psi(V)$.

The converse is also true (see \cite[Proposition~2.6]{SplittingThmEulerlike}). We paraphrase it in the version given in~\cite[Theorem~5.1]{HigsonReza}.

\begin{Theorem}\label{thm:TNE}
Let $X\in \vX(M)$ be a Euler-like vector field for $N$. There exists a neighborhood~$V$ of the zero section of $\nu N$ and a tubular neighborhood embedding $\psi\colon V \to M$ such that the Euler vector field is mapped to $X|_{\psi(V)}$.
\end{Theorem}

\begin{Remark}
Given $X$, the germ of such a tubular neighborhood embedding is unique.
\end{Remark}

A geometric proof of this result is given in
\cite{HigsonReza}, see also \cite[Sections~2.2 and~2.3]{DefSpacesEulerLike}. We review a few key-points following the latter reference.
\begin{itemize}\itemsep=0pt
\item Given the submanifold $N$ of $M$, the associated deformation space
is the set
\[
\cD(M,N):=\nu N \sqcup \bigl(M\times \RR^{\times}\bigr),
\]
where $\RR^{\times}=\RR\setminus\{0\}$. This set is endowed with a certain manifold structure.

\item Every function $f\in I_N$ vanishing on $N$ induces a smooth function $\widetilde{f}$ on the deformation space $\cD(M,N)$, which agrees with $({\rm d}f)|_N$ on $\nu N$ and with
$t^{-1}f $ on $M\times \RR^{\times}$, {where $t$ is the standard coordinate on $\RR$.}

\item Every $Y\in \vX(M)$ tangent to $N$ induces a vector field $\cD(Y)$ on the deformation space, which reads $Y_{\rm lin}$ on $\nu N$ and $Y\times 0$ on $M\times \RR^{\times}$.

\item An Euler-like vector field $X$ gives rise not only to $\cD(X)$ but also to
another vector field~$W$ on the deformation space $\cD(M,N)$, which reads \smash{$\pd{t}+\frac{1}{t}X$} on
$M\times \RR^{\times}$.

If we assume that the Euler-like vector field $X$ is complete, then
the time $1$ flow $\phi_W^1$ of~$W$ {is defined on $\nu N$,} maps $\nu N$ to {$M\times \{1\}\cong M$}, and the resulting map is a tubular neighborhood embedding. It maps the Euler vector field to $X$, since one can check that~${[W,\cD(X)]=0}$.

{In the general case, since on $N\times \RR$ the vector field $W$ is given by \smash{$\pd{t}$}, there
the flow $\phi_W^t$ is certainly defined for all times $t\in \RR$. Hence the time $1$ flow $\phi_W^1$ is defined in a neighborhood of $N\times \RR$ in $\cD(M,N)$. In particular, it is defined on a neighborhood of $N$ in $\nu N$, yielding a tubular neighborhood embedding.}

{We also revisit briefly the construction of \cite[Section~5]{HigsonReza}, when
 $X$ is not complete. Let $C$ be an open subset} of $\nu N$ with compact closure. On $C$, the flow $\phi_W^s$ is defined for some~${s>0}$.
{A} tubular neighborhood embedding {defined on $sC$ is given by} \smash{$v\mapsto \phi_W^s\bigl(s^{-1} v\bigr)$}, which agrees with $\phi_W^1(v)$. The time $1$ flow $\phi_W^1$ does not depend on the choice of $s$, so it is defined on a~neighborhood of the zero section of $\nu N$, and is a tubular neighborhood embedding there.
\end{itemize}

\begin{Remark}[functoriality]\label{rem:funct}
The construction in Theorem~\ref{thm:TNE} is functorial, as explained in~\cite[Sections~2.2 and~2.4]{SplittingThmEulerlike}, and as we now explain for Euler-like vector fields as in Definition~\ref{def:tneeulerlike}.
Fix pairs $(M,N)$ and \smash{$\bigl(\widetilde{M},\widetilde{N}\bigr)$} consisting of a manifold and a submanifold. Fix a smooth map \smash{$F\colon \widetilde{M}\to M$} such that \smash{$F\bigl(\widetilde{N}\bigr)\subset N$}, then by taking derivatives there is an induced morphism \smash{$\nu F\colon \nu \widetilde{N}\to \nu N$} between the normal bundles. \big(It is an isomorphism on each fiber whenever $F$ is transverse to $N$ and \smash{$F^{-1}(N)=\widetilde{N}$}.\big)

Fix an Euler-like vector field \smash{$\widetilde{X}$} for \smash{$\widetilde{N}$}, and an Euler-like vector field $X$ for $N$, such that~\smash{$\widetilde{X}$} is $F$-related to $X$.
Let $\psi\colon V \to M$ be a tubular neighborhood embedding
associated to~$X$ as in Theorem~\ref{thm:TNE},
for some neighborhood $V$ of the zero section of $\nu N$.
Let \smash{$\widetilde{\psi}\colon \widetilde{V}\to \widetilde{M}$} a~tubular neighborhood embedding associated to \smash{$\widetilde{X}$}, such that \smash{$\widetilde{V}\subset (\nu F)^{-1}(V)$} (this can always be arranged shrinking the domain).
 Then \smash{$\widetilde{\psi}$}, $\psi$ intertwine $\nu F$ and $F$, i.e., this diagram commutes:
\begin{gather}\label{eq:diagcomm}
\begin{split}
& \xymatrix{
\widetilde{V} \ar[r]^{\widetilde{\psi}}\ar[d]_{\nu F} & \widetilde{M}\ar[d]^{F}\\
V \ar[r]^{\psi} & M
			}
\end{split}
\end{gather}
\end{Remark}

We show that in certain cases a tubular neighborhood embedding can be ``lifted'', by using the above functoriality.
\begin{Lemma}\label{lem:ehresmann}
Let \smash{$F\colon \widetilde{M}\to M$} be surjective submersion, {$N\subset M$ a submanifold},
and ${\psi} \colon V \to M$ a tubular neighborhood embedding defined on $V\subset \nu N$.
Then there exists
a tubular neighborhood embedding \smash{$\widetilde{\psi}\colon \widetilde{V}\to \widetilde{M}$}, where \smash{$\widetilde{N}:=F^{-1}(N)$} and \smash{$\widetilde{V}\subset (\nu F)^{-1}(V)$}, such that
diagram \eqref{eq:diagcomm} commutes.
\end{Lemma}
\begin{proof}
Denote by ${X}$ the Euler-like vector field {on $\psi(V)$} corresponding to $\psi$.
 Take an Ehresmann connection for the submersion $F$. Then the lift $\widetilde{X}$ is
a vector field\footnote{Even if $X$ is complete, the lift $\widetilde{X}$ is not necessarily complete, as one sees considering the example $\widetilde{M}=\bigl\{(x,y)\in \RR^2: xy\in (-1,1)\bigr\}$, $M=\RR$, $F$ given by the first projection, $N=\{0\}$. This is {also} a reason why we do not require completeness in Definition \ref{def:tneeulerlike}\,(ii).}
 on $F^{-1}(\psi(V))$ which is $F$-related to $X$ and {satisfies} \smash{$\widetilde{X}|_{\widetilde{N}}=0$}.

We show that \smash{$\widetilde{X}_{\rm lin}$} is the Euler vector field of $\nu \widetilde{N}$.
For all \smash{$\widetilde{g}\in C^{\infty}\bigl(\widetilde{N}\bigr)$}, we have $\smash{\widetilde{X}_{\rm lin}(\widetilde{g})}= \smash{\bigl(\widetilde{X}|_{\widetilde{N}}\bigr)(\widetilde{g})}=0$.
For all $f\in I_N$, we have
\[
\widetilde{X}_{\rm lin}[F^*f]=\bigl[\widetilde{X}(F^*f)\bigr]=[F^*(X(f))]=(\nu F)^*[X(f)]=(\nu F)^*[f]=[F^*f],
\]
where we used that $[F^*h]=(\nu F)^*[h]$
 for all $h\in I_N$ in the third and last equality, and $[X(f)]=X_{\rm lin}[f]=E[f]=[f]$ in the fourth one.
Remark~\ref{rem:charEuler} then implies that $\widetilde{X}_{\rm lin}$ is the Euler vector field.

Hence \smash{$\widetilde{X}$} is an Euler-like vector field for \smash{$\widetilde{N}$}. Let \smash{$\widetilde{\psi}$} be a tubular neighborhood embedding associated to it.
To conclude the proof, we use the functoriality of the tubular neighborhood embedding construction recalled in Remark~\ref{rem:funct}.
\end{proof}

\subsection{A linearization criterion}

We present an equivalent characterization of the linearizability of a singular foliation around an embedded leaf.
When the leaf is just a point, the statement already appears in \cite[Section~8]{CerveauSing}.
\begin{Proposition}\label{prop:Cerveau}\samepage
Let $(M,\cF)$ be a foliated manifold, $L$ an embedded leaf. Then $\cF$ is linearizable around $L$ if and only if
\begin{itemize}\itemsep=0pt
\item[$(a)$] for some point $p\in L$, the map of equation~\eqref{eq:notinj} is injective,
\[\lambda_p \colon \ \g_p \to \vX_{\rm lin}(\nu_p L)=\operatorname{End}(\nu_p L),\qquad Y \bmod I_p\cF \mapsto Y_{\rm lin}|_{\nu_p L},
\]
\item[$(b)$] in a neighborhood of $L$ there is Euler-like vector field~$X$ for~$L$ such that $[X,
\cF]\subset \cF$.
\end{itemize}
\end{Proposition}

\begin{Remark}\label{rem:allpoints}
If condition (a) holds at a point of $L$, then it holds at all points of $L$, due to the involutivity of $\cF$.
\end{Remark}

\begin{Remark}\label{rem:refhom}
Proposition~\ref{prop:Cerveau} is related to \cite[Theorem~1.32]{NbhdSingLeaf}, which holds when the singular foliation $\cF$ is
locally real analytic, i.e., when nearby every point it admits real analytic generators. Indeed, the second item of \cite[Theorem~1.32]{NbhdSingLeaf} states that when condition~(b) in the above proposition is satisfied, near $L$ the
singular foliation~$\cF$ is generated by vector fields which are homogeneous (with respect to the
Euler-like vector field~$X$). Proposition~\ref{prop:Cerveau} states that when in addition condition a) is satisfied, $\cF$ is generated by vector fields which are homogeneous of degree~$1$.
\end{Remark}

\begin{proof}
``$\Rightarrow$'' Let $\psi\colon V\to M$ be a tubular neighborhood embedding such that $\psi^{-1}\cF=\cF_{\rm lin}$, {where $V$ is a neighborhood of the zero section of $\nu L$}. For any $p\in L$, we obtain a submanifold $S_p:=\psi(V_p)$
transverse to $L$, where $V_p:=V\cap \nu_p L$. The restricted foliation
 \smash{$\cF|_{S_p}:=\iota^{-1}_{S_p}\cF$} is linearizable, by means of the restriction of $\psi$ to a map $V_p\to S_p$. Hence
 \[
 \cF|_{S_p}/I_p\bigl(\cF|_{S_p}\bigr)\to
\bigl(\cF|_{S_p}\bigr)_{\rm lin}/I_0\bigl(\cF|_{S_p}\bigr)_{\rm lin}, \qquad [Y]\mapsto [Y_{\rm lin}]
\]
is an isomorphism.

This isomorphism is just the map $\lambda_p$ in the statement, restricted to its image. Therefore, $\lambda_p$ is injective, yielding (a).
To see that the maps agree, use $\cF|_{S_p}/I_p\bigl(\cF|_{S_p}\bigr)=\g_p $
(by the splitting theorem for singular foliations \cite[Theorem~1.1]{CerveauSing} and \cite[Section~1.3]{AndrSk})
and $\bigl(\cF|_{S_p}\bigr)_{\rm lin}=(\cF_{\rm lin})|_{\nu_p L}$.

Denote $X:=\psi_*(E)$, an Euler-like vector field for $L$ defined on $\psi(V)$.
Since the pushforward of vector fields preserves the Lie bracket, the condition $[X,\cF]\subset \cF$ {on $\psi(V)$} is equivalent to $[E,\cF_{\rm lin}]\subset \cF_{\rm lin}$. The latter condition holds because $\cF_{\rm lin}$ is generated by
fiberwise linear vector fields, and the Euler vector field $E$ commutes with all fiberwise linear vector fields (as we saw in Remark~\ref{rem:charEuler}). This yields (b).

``$\Leftarrow$'' We show that a tubular neighborhood embedding $\psi\colon V\to M$, associated as in Theorem~\ref{thm:TNE} to the Euler-like vector field $X$,
linearizes the singular foliation $\cF$.
We proceed analogously\footnote{The main difference is that we linearize along a leaf, rather then along a transversal, and that the Euler-like vector field $X$ does not belong to $\cF$.}
 to the proof of \cite[Proposition~2.6]{DefSpacesEulerLike}.

 Denote by
$\cD(\cF)$ the singular foliation on $\cD(M,L)$ generated by $\{\cD(Y): Y\in \cF\}$.
It is really a singular foliation: involutivity follows from the relation
$[\cD(Y),\cD(Z)]=\cD[Y,Z]$, which clearly holds on the dense subset $M\times \RR^{\times}$ and by continuity on the whole deformation space. The fact that $\cD(\cF)$ is locally finitely generated follows from the fact that $\cF$ is.

Denote by $W$ the vector field that $X$ induces on $\cD(M,L)$, as recalled at the beginning of Section~\ref{subsec:tne}.
We claim that
\begin{equation}\label{eq:wDF}
[W,\cD(\cF)]\subset \cD(\cF).
\end{equation}
 This claim and Proposition~\ref{prop:OriAlfcompl} imply that the flow of $W$ preserves $\cD(\cF)$. The restriction of~$\cD(\cF)$ to $\nu L$ and $M\times \{s\}$ (for any $s>0$) yields $\cF_{\rm lin}$ and $\cF$ respectively.
On any open subset $C$ of~$\nu L$ with compact closure, the flow $\phi_W^s$ is defined for some $s>0$. The dilation $s^{-1} \cdot\colon \nu L\to \nu L$ preserves $\cF_{\rm lin}$, because the Euler vector field does (a consequence of Remark~\ref{rem:charEuler}). By the description of the tubular neighborhood embedding $\psi$ in Section~\ref{subsec:tne}, we infer that its restriction to $C$ maps $\cF_{\rm lin}$ to $\cF$, proving the proposition.

We are left with proving the claim \eqref{eq:wDF}. Fix $Y\in \cF$. We have
\begin{equation}\label{eq:t-1}
[W,\cD(Y)]=t^{-1}\cD[X,Y]
\end{equation}
 since this clearly holds on the dense subset $M\times \RR^{\times}$ of the deformation space.
By assumption~(b), we have $[X,Y] \in \cF$. We now show the stronger statement $[X,Y] \in I_L\cF$.

Consider the module map $\overline{\rm lin}$ and the corresponding vector bundle map $\lambda$, both introduced in Section~\ref{subsec:linfol}. Note that we have the following commutative diagram of vector bundles (actually, Lie algebroids) over $L$ with exact rows:
\[
\begin{tikzcd}
0 \arrow{r} &\ker(\rho_{A_L}) \arrow{d}\arrow{r}& A_L \arrow{r}{\rho_{A_L}}\arrow{d}{\lambda}&TL \arrow{r}\arrow{d}& 0 \\
0\arrow{r} & \operatorname{End}(\nu L) \arrow{r}& {{\rm CDO}(\nu L)} \arrow{r}{\sigma} & TL \arrow{r} & 0.
\end{tikzcd}
\]
The left vertical map is injective, by assumption (a): to see this, use $\ker(\rho_{A_L})|_p=\g_p$, together with the fact that {$\lambda$} has constant rank (see Remark~\ref{rem:flin}) {or alternatively with Remark~\ref{rem:allpoints}}. Since the right vertical map is the identity, it follows that the middle vertical map {$\lambda$} is injective, which implies that $\overline{\rm lin}$ in \eqref{eq:overlinelin} is injective.
Under $\overline{\rm lin}$, the class of $[X,Y]$ is mapped to zero,
as
$
[X,Y]_{\rm lin} = [X_{\rm lin},Y_{\rm lin}] = 0
$
(the Euler vector field $X_{\rm lin}$ commutes with all fiber-wise linear vector fields).
The injectivity of $\overline{\rm lin}$ now implies that $[X,Y]\in I_L\cF$.

So we can write $[X,Y]=\sum_i f_iY_i$ as a finite sum of functions $f_i\in I_L$ and elements $Y_i\in \cF$. Therefore, $\cD([X,Y])=\sum_i f_i\cD(Y_i)$ on $M\times \RR^{\times}$, and \[t^{-1}\cD([X,Y])=\sum_i \widetilde{f}_i\cD(Y_i),\]
where the functions $\widetilde{f}_i$ on the deformation space were introduced in Section~\ref{subsec:linfol}.
This vector field lies in $\cD(\cF)$, since the functions $\widetilde{f}_i$ are smooth.
Together with equation~\eqref{eq:t-1} this shows that $[W,\cD(Y)]\in \cD(\cF)$, implying the claim \eqref{eq:wDF}.
 \end{proof}

\subsection{Hausdorff Morita equivalence}

We recall some facts about Hausdorff
Morita equivalence for singular foliations, which was introduced in \cite{HME}.
\begin{Definition}[{\cite[Definition~2.1]{HME}}]
Two foliated manifolds $(M_1,\cF_1)$ and $(M_2,\cF_2)$
are {\it Hausdorff
Morita equivalent} if there exist a manifold $P$ and two \emph{surjective submersions with connected fibers} $\pi_1 \colon P\to M_1$ and $\pi_2 \colon P\to M_2$ such that the pullback foliations agree: $\pi_1^{-1}\cF_1=\pi_2^{-1}\cF_2$. In this case, we write $(M_1,\cF_1)\simeq_{ME} (M_2,\cF_2)$,
\[\begin{tikzcd}
 & P \arrow[dl,swap, "\pi_1"] \arrow[dr, "\pi_2"] &\\
 (M_1,\cF_1)& &(M_2,\cF_2).
\end{tikzcd}\]
\end{Definition}

\begin{Remark}
When two foliated manifolds are Hausdorff
Morita equivalent, they share the same ``global transverse geometry''. For instance, their leaf spaces are homeomorphic, in a~way that preserves the property of being an embedded leaf, and the codimension of leaves;
 the isotropy Lie algebras\footnote{As well as the isotropy Lie groups of the holonomy groupoids associated to the singular foliations.}
at corresponding leaves
are isomorphic; their normal representations are isomorphic \cite[Proposition~2.5 and Theorem~3.44]{HME}.

In particular, let $(P,\pi_1,\pi_2)$ be a Hausdorff
Morita equivalence between the foliated manifolds $(M_1,\cF_1)$ and $(M_2,\cF_2)$. Let $L_1\subset M_1$ be a leaf of $\cF_1$. Then the corresponding
leaf $L_2$ of $\cF_2$ is determined by the property $\pi_1^{-1}(L_1)=\pi_2^{-1}(L_2)$.
\end{Remark}

\begin{Remark}
Let $G_i$ be a source-connected, Hausdorff Lie groupoid over a manifold $M_i$, for~${i=1,2}$. Denote by $\cF_i$ the induced singular foliations on $M_i$ as in Section~\ref{sec:singfol}.

{\samepage Assume that $G_1$ and $G_2$ are Morita equivalent as Lie groupoids.
 Then
\begin{enumerate}\itemsep=0pt
\item[(a)] The singular foliations $\cF_1$ and $\cF_2$ are
Hausdorff Morita equivalent \cite[Proposition~2.29]{HME}.
\item[(b)] The Lie groupoid $G_1$ is linearizable around an embedded leaf $L_1 \subset M_1$ if and only if
the Lie groupoid $G_2$ is linearizable around the corresponding leaf $L_2$ in $M_2$ \cite[Proposition~3.7]{LinProperGroupoids}.
\end{enumerate}
 In} that case, it follows that both $\cF_1$ and $\cF_2$ are linearizable along these leaves; indeed, the singular foliation on $\nu L_i$ induced by the linearization of $G_i$ is $(\cF_i)_{\rm lin}$, $i=1,2$ \cite[Lemma 4.15]{AZ2}.
(The converse is not true: a singular foliation can be linearizable, even when it arises from a Lie groupoid which is not linearizable.)

These two facts {combine into a hint toward} the main statement of this note, namely Theorem~\ref{thm:main} in the next section.
\end{Remark}

\section{The main theorem}

The following is the main result of this note. This section is dedicated to its proof.

\begin{Theorem}\label{thm:main}
Let $(M_1,\cF_1)$ and $(M_2,\cF_2)$ be Hausdorff Morita equivalent singular foliations,
and $L_1\subset M_1$, $L_2\subset M_2$ corresponding embedded leaves.
Then $\cF_1$ is linearizable around $L_1$ if and only if $\cF_2$ is linearizable around $L_2$.
\end{Theorem}

The theorem follows immediately from the following proposition.
\begin{Proposition}\label{prop:pss}
 Let $\pi\colon P\to M$ be a surjective submersion with connected fibers. Let $\cF$ be a~singular foliation on $M$ and $L\subset M$ an embedded leaf.
 Then
 $\cF$ is linearizable around $L$ if and only if $\pi^{-1}\cF$ is linearizable around $\pi^{-1}(L)$.
\end{Proposition}

To prove Proposition~\ref{prop:pss}, we first need a lemma, for which {the connectedness of fibers is not used, and for which} it might be useful to keep in mind this diagram
\begin{equation*}
\xymatrix{
\nu \bigl(\pi^{-1}(L)\bigr) \ar[r]^{}\ar[d]_{\nu \pi} & \pi^{-1}(L)\; \ar[d]^{} \ar@{^{(}->}[r]& P\ar[d]^{\pi} \\
\nu L \ar[r]^{} & L\;\ar@{^{(}->}[r]& M.
			}
\end{equation*}

\begin{Lemma}\label{lem:piFlin}
Assume the set-up of Proposition~$\ref{prop:pss}$.
\begin{itemize}\itemsep=0pt
\item[$(i)$] If $Z\in \vX(P)$ is {tangent to $\pi^{-1}(L)$ and} $\pi$-projectable to $X\in \vX(M)$, then $Z_{\rm lin}$ is $\nu \pi$-projectable to $X_{\rm lin}$.
\item[$(ii)$] We have
$
\bigl(\pi^{-1}\cF\bigr)_{\rm lin}=(\nu \pi)^{-1}(\cF_{\rm lin})$.
\end{itemize}
\end{Lemma}

\begin{proof}
For item (i), we first have to check that
$ Z_{\rm lin}((\nu \pi)^*(g))=(\nu \pi)^*(X_{\rm lin}(g))$
for all $g\in C^{\infty}(L)$ (viewed as fiberwise constant functions on $\nu L$).
This holds by a direct computation, using that~$Z|_{\pi^{-1}(L)}$ is projectable to $X|_L$.
We also have to check the same equation
replacing $g$ by the fiberwise linear function $[f]$ for all $f\in I_L$ (here we adopt the notation of Section~\ref{subsec:linfol}). This is done similarly, using that $(\nu \pi)^*[f]=[\pi^*f]\in C_{\rm lin}\bigl(\nu(\pi^{-1}L)
\bigr)$, as in proof of Lemma~\ref{lem:ehresmann}.

Item (ii) is an equality of modules of vector fields on $\nu \bigl(\pi^{-1}(L)\bigr)$. Recall that{\samepage
\begin{itemize}\itemsep=0pt
\item $\bigl(\pi^{-1}\cF\bigr)_{\rm lin}$ is generated by $Z_{\rm lin}$ for vector fields $Z$ which $\pi$-project to elements of $\cF$.
\item $(\nu \pi)^{-1}(\cF_{\rm lin})$ is generated by vector fields $Y$ which $\nu \pi$-project to elements of $\cF_{\rm lin}$.
\end{itemize}
We see that item (i) immediately implies ``$\subset$'' in item (ii).}

For the inclusion ``$\supset$'', take a vector field $Y$ on $\nu\bigl(\pi^{-1}(L)\bigr)$ which $\nu \pi$-projects to an element of $\cF_{\rm lin}$; the latter is necessarily a finite sum $\sum_i h_i (X_i)_{\rm lin}$ for some $h_i\in C_c^{\infty}(\nu L)$ and
$X_i\in \cF$. Since $\pi$ is a surjective submersion, there are $Z_i\in \pi^{-1}\cF$ which $\pi$-project to the $X_i$; by item (i),
\[
S:=\sum_i ((\nu \pi)^*h_i) (Z_i)_{\rm lin}
\]
is $\nu \pi$-projectable to $\sum_i h_i (X_i)_{\rm lin}$. Therefore, the difference $Y-S$ is tangent to the fibers of $\nu \pi$.
As~$Y=(Y-S)+S$ and \smash{$S\in \widehat{(\pi^{-1}\cF)_{\rm lin}}$}, to finish the proof we just need to show that
\[
\Gamma_c(\ker((\nu \pi)_*))\subset \bigl(\pi^{-1}\cF\bigr)_{\rm lin}.
\]
Notice that $\ker((\nu \pi)_*)$ is a regular distribution, since $\nu \pi$ is a surjective submersion. It suffices to show that locally there are frames for
$\ker((\nu \pi)_*)$ consisting of elements of $\bigl(\pi^{-1}\cF\bigr)_{\rm lin}$. Such frames can be constructed as follows. Take a local frame $\{W_i\}$ {of $\ker(\pi_*)$ on $P$ and}
notice that these vector fields lie in $\pi^{-1}\cF$. We know that the $(W_i)_{\rm lin}$ lie in $\ker((\nu \pi)_*)$, by item (i), and they form a local frame for the latter.
\end{proof}

\begin{Remark}\label{rem:reffunct}
Most of Lemma~\ref{lem:piFlin} can also be proven using the deformation spaces recalled in Section~\ref{subsec:tne} and their functoriality \cite[Section~2.2]{SplittingThmEulerlike}: there is a natural smooth map
\[
\cD(\pi)\colon\ \cD(P,\pi^{-1}(L))\to \cD(M,L),
\]
which on the fiber over zero restricts to $\nu \pi\colon \nu \bigl(\pi^{-1}(L)\bigr)\to \pi^{-1}(L)$.
For item (i), the fact that $Z$ is $\pi$-projectable to $X$ immediately implies that~$\cD(Z)$ is $\cD(\pi)$-projectable to~$\cD(X)$ on an~open dense subset of $\cD\bigl(P,\pi^{-1}(L)\bigr)$, hence everywhere, in particular also on the fiber over zero. As~$\cD(\pi)$ is a surjective submersion, this
implies \smash{$\cD\bigl(\pi^{-1} \cF\bigr)\subset (\cD(\pi))^{-1} (\cD(\cF))$}, and restricting to the fiber over zero we obtain the inclusion ``$\subset$'' in item~(ii).
 \end{Remark}

{The following remark follows from Lemma~\ref{lem:piFlin}, and obviously it is consistent with Theorem~\ref{thm:main}.}

\begin{Remark}[Hausdorff Morita equivalence of the linearized foliations]
Let $(M_1,\cF_1)$ and $(M_2,\cF_2)$ be Hausdorff Morita equivalent singular foliations,
and $L_1$, $L_2$ be corresponding embedded leaves.
Then $(\nu L_1, (\cF_1)_{\rm lin})$ and $(\nu L_2, (\cF_2 )_{\rm lin})$ are Hausdorff Morita equivalent.

This follows immediately from Lemma~\ref{lem:piFlin}\,(ii): if $P$ is a manifold with surjective submersions with connected fibers $\pi_i$ to $M_i$ such that \smash{$\pi_1^{-1}\cF_1=\pi_2^{-1}\cF_2$},
then the manifold $\nu \tilde{L}$ and the maps $\nu \pi_i$ provide the desired Hausdorff Morita equivalence \big(in particular, $(\nu \pi_1)^{-1} (\cF_1)_{\rm lin}=(\nu \pi_2)^{-1} (\cF_2)_{\rm lin}$\big). Here \smash{$\tilde{L}:=\pi_1^{-1}(L_1)=\pi_2^{-1}(L_2)$}.
\end{Remark}

\begin{proof}[Proof of ``\boldmath{$\Rightarrow$}'' in Proposition~\ref{prop:pss}]
Suppose $\cF$ is linearizable around $L$.
Take a tubular neighborhood embedding $\psi\colon V \to M$ so that $\psi^{-1}\cF=\cF_{\rm lin}$, {where $V$ is a neighborhood of the zero section of $\nu L$}. One can ``lift'' $\psi$ to
a tubular neighborhood embedding $\widetilde{\psi}\colon \widetilde{V}
\to P$ around $\pi^{-1}(L)$ such that diagram \eqref{eq:diagcomm} commutes, by Lemma~\ref{lem:ehresmann}.

The commutativity of the diagram, {together with the functoriality of the pullback \cite[Proposition~1.11\,(b)]{AndrSk},} implies that $\widetilde{\psi}$ identifies $\pi^{-1}\cF$ with $(\nu \pi)^{-1}(\cF_{\rm lin})$. The latter equals $\bigl(\pi^{-1}\cF\bigr)_{\rm lin}$, by Lemma~\ref{lem:piFlin}\,(ii).

An alternative argument is as follows. Let $\widetilde{X}$ be the push-forward by
$\widetilde{\psi}$ of the Euler vector field of $\nu \bigl(\pi^{-1}(L)\bigr)$. One can show\footnote{If $\widetilde{Y}\in \vX(P)$ is $\pi$-projectable to $Y\in \cF$, then $\bigl[\widetilde{X}, \widetilde{Y}\bigr]$ $\pi$-projects to $[X,Y]\in \cF$, by the naturality of the Lie bracket. This implies the statement, thanks to Definition \ref{rem:pulblackfol} and the Leibniz rule.}
 that \smash{$\bigl[\widetilde{X},\pi^{-1}\cF\bigr]\subset \pi^{-1}\cF$}, using the fact that the Euler-like vector field $X$ associated to $\psi$ satisfies $[X,\cF]\subset \cF$ by the proof of Proposition~\ref{prop:Cerveau}. Then one can apply Proposition~\ref{prop:Cerveau} to $\bigl(P, \pi^{-1}\cF\bigr)$, to conclude again that $\pi^{-1}\cF$ is linearizable around~$\pi^{-1}(L)$. (Notice that hypothesis a) of Proposition~\ref{prop:Cerveau} is satisfied, since slices in $P$ transverse to~$\pi^{-1}(L)$ map diffeomorphically, as foliated manifolds, to slices in $M$ transverse to $L$).
\end{proof}

\begin{proof}[Proof of ``\boldmath{$\Leftarrow$}'' in Proposition~\ref{prop:pss}]

Suppose $\pi^{-1}\cF$ is linearizable around $\pi^{-1}(L)$.
We denote \smash{$\widetilde{\cF}:=\pi^{-1}\cF$} and \smash{$\widetilde{L}:=\pi^{-1}(L)$}. By assumption, there exists a tubular neighborhood embedding \smash{$\widetilde{\psi}\colon \widetilde{V} \to P$} such that \smash{$\bigl(\widetilde{\psi}\bigr)^{-1}\widetilde{\cF}=\bigl(\widetilde{\cF}\bigr)_{\rm lin}$},
where $\widetilde{V}$ is a neighborhood of the zero section of $\nu \widetilde{L}$. Denote by $\widetilde{X}$ the Euler-like vector field on the image of~$\widetilde{\psi}$, obtained as the push-forward by
$\widetilde{\psi}$ of the Euler vector field of $\nu \widetilde{L}$.
By Proposition~\ref{prop:Cerveau} (and its proof),
\begin{itemize}\itemsep=0pt
\item [(A)] for some point $p\in \w{L}$, the map
\[
\lambda_p \colon  \ \g_p \to \vX_{\rm lin}\bigl(\nu_p \w{L}\bigr),\qquad Y \bmod I_p\w\cF \mapsto Y_{\rm lin}|_{\nu_p \w L}
\]
is injective, where $\g_p$ is the isotropy Lie algebra of $\w \cF$ at $p$,
\item [(B)] $\bigl[\w{X},\w{\cF}\bigr]\subset \w{\cF}$.
\end{itemize}

{We want to show that the foliated manifold $(M,\cF)$ satisfies the two conditions of Proposition~\ref{prop:Cerveau}.}
For the foliated manifold $(M,\cF)$, condition (a) of Proposition~\ref{prop:Cerveau} is satisfied at the point $\pi(p)\in L$, for the same reason outlined in the previous proof: for any slice \smash{$\w{S}$} through $p$ transverse to \smash{$ \widetilde{L}$}, the restricted foliations \smash{$\bigl(\w{S},\iota_{\w{S}}^{-1} \w{\cF}\bigr)$} and \smash{$\bigl({S},\iota_{S}^{-1} {\cF}\bigr)$} are isomorphic via~$\pi$, where~\smash{${S:=\pi\bigl(\w{S}\bigr)}$}.

To check that $(M,\cF)$ satisfies condition (b) of Proposition~\ref{prop:Cerveau},
in the remainder of this proof we will construct an Euler-like vector field $X$ for $L$ such that $[X,\cF]\subset \cF$. Proposition~\ref{prop:Cerveau} will then assure the linearizability of $\cF$ around $L$.

\emph{Construction of X:} Notice that $\w{X}$ might not be $\pi$-projectable. To circumvent this difficulty we proceed as follows. Since \smash{$\pi|_{\w{L}} \colon \w{L} \to L$} is a surjective submersion, we can choose an open cover~$\{L_j\}$ of $L$ admitting sections \smash{$s_j\colon L_j\to \w{L}$}. Then
\[
\w{U}_j:=\w{\psi}\bigl({\widetilde{V}}|_{{\rm im}(s_j)}\bigr)
\]
is a submanifold of $P$
with \smash{$\dim\bigl(\w{U}_j\bigr)=\dim(M)$}, since $\w{L}$ and $L$ have the same codimension. Notice that at any point $q\in {\rm im}(s_j)$, the derivative of $\pi$ maps \smash{$T_q \w{U}_j=\w{\psi}_*\bigl(\widetilde{V}_q\bigr)\oplus T_q ({\rm im}(s_j))$} injectively (and thus isomorphically) into $T_{\pi(q)}M$. Therefore, \smash{$\w{U}_j$} is transverse to fibers of $\pi$ along~${\rm im}(s_j)$.
Shrinking \smash{$\w{U}_j$} if necessary, by the inverse function theorem, we have that \smash{$U_j:=\pi\bigl(\w{U}_j\bigr)$} is an~open neighborhood of $L_j$ in $M$ and
\[
\pi|_{\w{U}_j}\colon\ \w{U}_j\to U_j
\]
is a diffeomorphism.
Further this diffeomorphism matches the singular foliations \smash{$\iota_{\w{U}_j}^{-1}\w{\cF}$} and \smash{$\iota_{U_j}^{-1}{\cF}$}, by the functoriality of the pullback.

	\begin{figure}[!ht]
\centering
			\begin{tikzpicture}[scale=1]
			\draw[fill=gray!15!white] (-5.5,-1) -- (-3.5,1) -- (4.5,1) -- (2.5,-1) -- cycle; 
			\node[below,right] at (3.7,0) {$M$};
			
\draw[blue] (-1.5,-1) -- (0.5,1);		
\node[blue] at (0,0) {$L$};	
			
			\draw[-] (-1.5,0.5) -- (-1.5,5.5)--(0.5,7.5) -- (0.5,2.5) --cycle;
				\node[above] at (0.8,6) {$\widetilde{L}$};
				
 \draw[thick,cyan] (-1.5,4.5) .. controls (-1,5.5) and (-0.5,4) .. (-0.25,5);
 \node[right,cyan] at (-1.5,5.2) {${\rm im}(s_i)$};
 \draw[thick,purple] (-0.75,3) .. controls (-0.5,3.5) and (0,2.5) .. (0.5,4);
 \node[purple] at (-0.1,2.8) {${\rm im}(s_j)$};
			\end{tikzpicture}
		\caption{The images of two local sections $s_i$ and $s_j$ of $\pi|_{\w{L}} \colon \w{L} \to L$.}
	\end{figure}

Define $X_j\in \vX(U_j)$
so that it corresponds to \smash{$\w{X}|_{\w{U}_j}$} under the diffeomorphism \smash{$\pi|_{\w{U}_j}$}; this is~possible since \smash{$\w{X}$} is tangent to \smash{$\w{U}_j$}, {by the construction of the latter}. Using item (B) above, for all~$j$ we have that
\begin{itemize}\itemsep=0pt
\item $X_j$ is a Euler-like vector field for $L_j$,
\item $\bigl[X_j,\iota^{-1}_{U_j}\cF\bigr]\subset \iota^{-1}_{U_j}\cF$.
\end{itemize}
Let $\{\rho_j\}$ be a {(smooth and locally finite)} partition of unity on the open subset $\cup_j U_j\subset M$, subordinate to the covering $\{U_j\}$; {we may assume that each $\rho_j $ has compact support}.
 We define
\[X:=\sum_j\rho_jX_j.\]
This is an Euler-like vector field for $L$, as remarked in \cite[Section~2.4]{DefSpacesEulerLike}: one checks that \smash{$\bigl(\sum_j\rho_jX_j\bigr)_{\rm lin}=\sum_j(\rho_j)|_L (X_j)_{\rm lin}$} {is the Euler vector field}.

\emph{Checking the condition $[X,\cF]\subset \cF$:}
The failure of $\w{X}$ to be $\pi$-projectable is controlled, as~follows.
All vertical vector fields $V\in \Gamma_c(\ker \pi_*)$ lie in $\w{\cF}$, so by Proposition~\ref{prop:uptoF} -- which we can apply since
{$\w{X}$ is an infinitesimal symmetry of $\w{\cF}$}
by item (B) above --
the time-1 flow \smash{$\phi^V_1$} satisfies%
\begin{equation}\label{eq:difference}
\bigl(\phi^V_1\bigr)_*\w{X}-\w{X}\in \widehat{\w{\cF}}.
\end{equation}

 We claim that the difference ${X_i|_{U_i\cap U_j}-X_j|_{U_i\cap U_j}}$ lies in \smash{$\widehat{\cF}|_{U_i\cap U_j}$}, for all $i$, $j$ such that ${U_i\cap U_j\neq \varnothing}$. Indeed, since the $\pi$-fibers are connected, for every sufficiently small open subset~\smash{$\widetilde{W}\subset \w{U}_i\cap\bigl(\pi^{-1}( U_j)\bigr)$} there is a vertical vector field $V\in \Gamma_c(\ker \pi_*)$ whose time-1 flow maps \smash{$\widetilde{W}$} into
 \smash{$\w{U}_j\cap\bigl(\pi^{-1}(U_i)\bigr)$}. Let us assume for simplicity that \smash{$\widetilde{W}= \w{U}_i\cap\bigl(\pi^{-1}( U_j)\bigr)$}; then this diagram of diffeomorphisms commutes:
\begin{equation*}
\xymatrix{
 \w{U}_i\cap\bigl(\pi^{-1}( U_j)\bigr) \ar[rr]^{\phi^V_1 }\ar[dr]_{\pi} &&\ar[dl]_{\pi} \w{U}_j\cap\bigl(\pi^{-1}(U_i )\bigr)\\
& U_i\cap U_j& 			}
\end{equation*}
and the claim then follows\footnote{{Without the simplifying assumption, the above argument first yields a local version of the claim, namely: there is an open cover
$\{W_{\alpha}\}$
 of $U_i\cap U_j$
 satisfying
\smash{$Y|_{W_{\alpha}}\in \widehat{\cF}|_{W_{\alpha}}$} for all $\alpha$, where $Y:={X_i|_{U_i\cap U_j}-X_j|_{U_i\cap U_j}}$. Let $f\in C^{\infty}_c(U_i\cap U_j)$, and take a finite subcover of its support. Take a subordinate partition of unity $\{\rho_{\alpha}\}$. As~${\rho_{\alpha}f\in
C^{\infty}_c(W_{\alpha})}$, we see that
 $fY=\sum_{\alpha} \rho_{\alpha}f Y$ is a finite sum of elements of $\cF$, hence lies in $\cF$. This shows that~\smash{${Y\in \widehat{\cF}|_{U_i\cap U_j}}$}, as claimed.}
}
 from \eqref{eq:difference}.

Let $p\in L$, fix an index $j_0$ such that $p\in U_{j_0}$. {Consider
 the neighborhood $U^p:=\cap U_j$ of $p$, where the intersection is over all indices $j$ satisfying $p\in U_j$.}
Then on this neighborhood we have
\[X=\sum_j \rho_j X_j=
 X_{j_0}+
 \sum_j \rho_j (X_j- X_{j_0}),\]
 where the last sum lies in \smash{$\widehat{\cF}|_{U^p}$} by the above claim.
 Since \smash{$\bigl[X_{j_0},\iota^{-1}_{U_{j_0}}\cF\bigr]\subset \iota^{-1}_{U_{j_0}}\cF$}, we conclude that $[X,\cF]\subset \cF$ on $U^p$. {Since the Lie bracket is local, we therefore have $[X,\cF]\subset \cF$ in a~neighborhood of $L$ in $M$.}
 \end{proof}

\section{Examples and an application}

We start presenting examples of Hausdorff Morita equivalent and of linearizable singular foliations. Then, using Theorem~\ref{thm:main}, in Section~\ref{subsec:restr}--\ref{subsec:restrHME} we obtain a condition under which the linearizability of a singular foliation $\cF$ can be reduced to the one of a simpler foliation, namely its restriction to a slice.

\subsection{Examples}

Examples of Hausdorff Morita equivalent singular foliations are given in \cite[Sections~2.3--2.5]{HME}. Here we just mention one, namely \cite[Corollary~2.17]{HME}.

\begin{Example}
Let two connected Lie groups $G_1$, $G_2$ act freely and properly on a
manifold $P$ with commuting actions. Then the following singular foliations are Hausdorff Morita equivalent:
\begin{itemize}\itemsep=0pt
\item the singular foliation on $P/G_1$ given by the induced $G_2$ action,
\item the singular foliation on $P/G_2$ given by the induced $G_1$ action.
\end{itemize}
\end{Example}

We now present some examples of linearizable foliations.

\begin{Example}[Log tangent bundle]
Let $L$ be a codimension one submanifold of $M$. Let~$\cF_{\rm log}$ consists of all compactly supported vector fields\footnote{These are the
 compactly supported sections of the log tangent bundle associated to $L$.} on $M$ which are tangent to $L$.
Then~$\cF_{\rm log}$ is linearizable around $L$. Indeed, take any tubular neighborhood embedding to identify a~neighborhood of $L$ in $M$ with a neighborhood of $L$ in the line bundle $\nu L$. The log-foliation of vector fields on $\nu L$ tangent to the zero section is linear, as it is generated by all linear vector fields on~$\nu _L$.
\end{Example}

\begin{Example}[elliptic tangent bundle]
Similarly, let $f\colon M\to \RR$ be a Morse--Bott function, i.e., the critical set $L$ of $f$ is a submanifold such that the normal Hessian of $f$ is non-degenerate. Assume that $L$ has codimension $2$, and that $f|_L=0$. Let $\cF_{\rm ell}$ consist of all compactly supported vector fields\footnote{These are the compactly supported sections of the elliptic tangent bundle associated to $f$.} which preserve the ideal generated by~$f$. Then $\cF_{\rm ell}$ is linearizable around $L$.
Indeed, by the Morse--Bott lemma as in \cite[Section~4.2]{ELisotropic}, there is a
tubular neighborhood embedding that identifies $f$ in a neighborhood of $L\subset M$ with the normal Hessian of $f$ (viewed as a degree~$2$ polynomial) in a neighborhood of $L \subset \nu L$. The elliptic foliation associated to this polynomial is linear. This can be seen trivializing $\nu L$ by means of a local orthonormal frame, and using the fact that the elliptic foliation on $\RR^2$ associated to the function $x^2+y^2$ is generated by the Euler vector field and the rotation vector field, see, e.g., \cite[Section~1]{WitteCohomElliptic}.
\end{Example}

\begin{Remark}
The following simple remark is used in the next example.
Let $V\to L$ be a vector bundle and $\cF_V$ a singular foliation on $V$ generated by linear vector fields.
Let $G$ be a discrete group acting freely and properly on the vector bundle $V\to L$ by vector bundle automorphisms, preserving the singular foliation $\cF_V$.
Then $V/G$ is a vector bundle over $L/G$, and the projection $\pi \colon V\to V/G$ is a vector bundle morphism and a covering map. The singular foliation $\cF_V$ induces a singular foliation $\cF$ on $V/G$, which is also generated by linear vector fields.
\end{Remark}

\begin{Example}[mapping tori]
Consider a linear singular foliation $\cF_S$ on $\RR^n$ (i.e., one induced by the action of a Lie subalgebra of $\mathfrak{gl}(n,\RR)$), the full foliation on $\RR$, and their product foliation~$\cF_V$ on~${V=\RR^n\times \RR}$. For any $\phi\in {\rm GL}(n,\RR)$ preserving $\cF_S$, the group $\ZZ$ acts on $V$, with $1\in \ZZ$ acting by $(p,t)\mapsto (\phi(p),t+1)$. On the mapping torus
\[
V/\ZZ= (\RR^n\times [0,1])/{\sim}, \qquad\text{where }(p,0)\sim (\phi(p),1),
\]
we obtain a foliation $\cF$, again generated by linear vector fields. Here are some low-dimensional examples.
\begin{enumerate}\itemsep=0pt
\item[(i)] If $\cF_S$ is the foliation by points of $\RR$, and $\phi=-{\rm Id}$, we obtain the M\"obius strip $V/\ZZ$ with a~foliation by circles. Most circles wind twice around the M\"obius strip, while the ``middle circle'' winds only once and has holonomy group $\ZZ_2$.

\item[(ii)] If $\cF_S$ is the foliation generated by \smash{$y\pd{y}$} on $\RR$, on the M\"obius strip, we obtain a foliation with precisely two leaves: a circle, and a 2-dimensional open leaf.

\item[(iii)] If $\cF_S$ is the foliation generated by the rotation vector field \smash{$x\pd{y}-y\pd{x}$} of $\RR^2$, and $\phi(x,y)=(x,-y)$ on $V/\ZZ$, we obtain a singular foliation whose leaves are a circle and many copies of the Klein bottle.
\end{enumerate}
\end{Example}

\begin{Remark}\label{rem:mappingtori}
Notice that in both in (i) and (ii) above, the foliation on the M\"obius strip and the foliation restricted to a slice $\RR$ are not Hausdorff Morita equivalent. For instance, in (ii), $\cF$~has two leaves but the restricted foliation has three.
\end{Remark}

\subsection{Restrictions of singular foliations to slices}\label{subsec:restr}

Let $(M,\cF)$ be a foliated manifold, $L$ an embedded leaf.

\begin{Definition}
A \textit{slice} at a point $p\in L$ is a submanifold $N$ of $M$ such that $N\cap L=\{p\}$, we have a direct sum $T_pN\oplus T_pL=T_pM$, and $N$ is transverse to the leaves of $\cF$ it meets.
\end{Definition}
In that case, \smash{$\cF|_N:=\iota^{-1}_N\cF$} is a singular foliation on $N$, with the point $\{p\}$ being a leaf. If~we take another slice $N'$ through $p$, upon shrinking the slices we obtain isomorphic foliated manifolds \cite[Appendix A.2]{AZ2}.

\begin{Lemma}\label{lem:nuplin}
If $\cF$ is linearizable along $L$, then $ \cF|_N$ is linearizable at $p$.
\end{Lemma}
\begin{proof}
Take a tubular neighborhood embedding $\phi \colon V\to M$ linearizing $\cF$, where $V$ is open in $\nu L$. Restrict $\phi$ to a tubular neighborhood embedding from $(V_p, \cF_{\rm lin}|_{V_p})$ to $(N', \cF|_{N'}) $, where $V_p:=V\cap \nu_pL$ and $N':=\phi(V_p)$. Then use $(\cF|_{N'})_{\rm lin}=\cF_{\rm lin}|_{V_p}$, which holds by the splitting theorem for singular foliations \cite[Theorem~1.1]{CerveauSing} and \cite[Section~1.3]{AndrSk}.
\end{proof}

The converse of Lemma~\ref{lem:nuplin} does not hold. That is, given
 a singular foliation $\cF$ and a slice~$N$ to an embedded leaf $L$ at $p$:
if $\cF|_{N}$ is linearizable at $p$, then in general
it does not follow that $\cF$ is linearizable along $L$.
This is well-known and already apparent in the case of regular foliations, as we recall in this example.

\begin{Example}
On the cylinder $S^1\times \RR$ with coordinates $\theta$ and $y$, consider the (involutive) distribution $D$ spanned by the vector field \smash{$\pd{\theta}+g(y)\pd{y}$},
where $g\colon \RR\to \RR$ is a smooth function which is identically zero for $y\le 0$ and which is $>0$ for $y>0$. It gives rise to a foliation $\cF=\Gamma_c(D)$, with $L=S^1\times \{0\}$ as a leaf. Take the slice $N=\{p\}\times \RR$ through a point $p\in S^1$; the restricted foliation is the zero foliation, which is certainly linearizable. But $\cF$ itself is not linearizable along $L$: the linearized foliation $\cF_{\rm lin}$ on $\nu L= S^1\times \RR$ is spanned by \smash{$\pd{\theta}$}, since $g'(0)=0$, and thus has no holonomy around the zero section, while $\cF$ does have holonomy around $L$, since $g$ is positive on $\{y>0\}$.

We mention that by a version of the Reeb stability theorem for regular foliations,
if a regular foliation has a (not necessarily compact) embedded leaf $L$ with finite holonomy group, then the foliation is linearizable around $L$, see, e.g., \cite[Theorem~3.2.1]{MarcutThesis}. If $L$ is \emph{compact}, the classical Reeb stability theorem assures that $L$ is linearizable by a tubular neighborhood embedding
with image a \emph{saturated} neighborhood of $L$, i.e., one given by a union of leaves.
\end{Example}

 \subsection{Hausdorff Morita equivalence and restrictions to slices}\label{subsec:restrHME}

We determine a condition under which the converse of Lemma~\ref{lem:nuplin} holds, using Theorem~\ref{thm:main}.
This is desirable because in that case the linearizability of $\cF$ becomes equivalent to the linearizability of its restriction to a slice; for the latter -- which vanishes at a point -- some linearization criteria are given by Cerveau in \cite{CerveauSing}.

In general, in a neighborhood of a leaf, a singular foliation is not
Hausdorff Morita equivalent to its restriction to a slice, as we saw in Remark~\ref{rem:mappingtori}. Below we determine conditions under which they are.

\begin{Remark}\label{rem:Hppconn}
 If in a neighborhood of a leaf $L$, a singular foliation $\cF$ is
Hausdorff Morita equivalent to its restriction \smash{$\iota^{-1}_N\cF$} to a slice $N$, then the holonomy group of $\cF$ at the point $\{p\}=N\cap L$ is \emph{connected}.

Indeed, $\cF$ and \smash{$\iota^{-1}_N\cF$} being Hausdorff Morita equivalent implies that their holonomy groups at $p$ are isomorphic as Lie groups \cite[Theorem~3.44\,(i)]{HME}.
The holonomy group of \smash{$\iota_N^{-1}\cF$} at $p$ is connected because it is the whole source fiber of the holonomy groupoid \cite{AndrSk} of \smash{$\iota_N^{-1}\cF$}, since the latter has $\{p\}$ as a leaf, and holonomy groupoids are always source-connected.
\end{Remark}

Recall that given a singular foliation on $M$ and a slice $N$, the saturation of $N$ is the union of all leaves intersecting $N$; it is an open subset of $M$.

 \begin{Proposition}\label{prop:MEslice}
 Let $\cF$ be a singular foliation arising $($in the sense of Section~$\ref{subsec:sf})$ from a~Hausdorff, source connected, Lie groupoid $\cG$.
Let $L$ be an embedded leaf, $N$ a slice. Suppose the restricted Lie groupoid $\cG_N:=\bt^{-1}(N)\cap \bs^{-1}(N)$ is source connected.
Then $\bigl(N,\iota_N^{-1}\cF\bigr)$ and $\bigl(U,\iota_U^{-1}\cF\bigr)$ are Hausdorff Morita equivalent singular foliations, where~$U$ is the saturation of $N$.
\end{Proposition}
\begin{proof}
Notice that $\widetilde{N}:=\bs^{-1}(N)$ is a smooth Hausdorff manifold. Consider the two maps
\[\begin{tikzcd}
 & \widetilde{N} \arrow[dl,swap, "\bt|_{\widetilde{N}}"] \arrow[dr, "\bs|_{\widetilde{N}}"] &\\
U& &N,
\end{tikzcd}\]
where $U=\bt\bigl(\widetilde{N}\bigr)$ agrees with the saturation of $N$.

Clearly, \smash{$\bs|_{\widetilde{N}}$} is a surjective submersion with connected fibers.

The map \smash{$\bt|_{\widetilde{N}}\colon \widetilde{N} \to U$} is a submersion, as a consequence of the fact that $N$ is transverse to the leaves it meets
(see, e.g., \cite[Example~4.2.1 and Proposition~5.2.4]{MatiasME}).

We show that \smash{$\bt|_{\widetilde{N}}$} has connected fibers. Let $u\in U$. The corresponding fiber of \smash{$\bt|_{\widetilde{N}}$} is \smash{$ \bt^{-1}(u)\cap \widetilde{N}$}. There exists $p\in N$ lying in the same leaf
 as $u$; take any $g\in \cG$ with $\bt(g)=p$ and $\bs(g)=u$. Then the left composition with $g$ gives a diffeomorphism
\smash{$ \bt^{-1}(u)\cap \widetilde{N}\to \bt^{-1}(p)\cap \widetilde{N}$}.
The latter is the target fiber of $\cG_N$ at $p$, so it is connected by assumption.

Finally, we show that
\[
\bigl(\bt|_{\widetilde{N}}\bigr)^{-1}\cF=\bigl(\bs|_{\widetilde{N}}\bigr)^{-1}\bigl(\iota_N^{-1}\cF\bigr).
\]
Since $\cG$ is a Lie groupoid giving rise to $\cF$, it is also a bisubmersion \cite[Proposition~2.2]{AndrSk}, i.e., $\bt^{-1}\cF=\bs^{-1}\cF$.
By the functoriality of pullback, we have
\[
\bigl(\bs|_{\widetilde{N}}\bigr)^{-1}\iota_N^{-1}\cF=\iota_{\widetilde{N}}^{-1} \bigl(\bs^{-1} \cF\bigr)=\iota_{\widetilde{N}}^{-1} \bigl(\bt^{-1} \cF\bigr)=\bigl(\bt\circ \iota_{\widetilde{N}}\bigr)^{-1} \cF.
\tag*{\qed}
\]
\renewcommand{\qed}{}
\end{proof}

\begin{Remark}
Proposition~\ref{prop:MEslice} can be also deduced as follows.
The inclusion $\cG_N\to \cG_U$ is a~weak equivalence of Lie groupoids, yielding that these two Lie groupoids are Morita equivalent. Hence, when both are source connected (and $\cG_U$ automatically is), their singular foliations are Hausdorff Morita equivalent, by \cite[Proposition~2.29]{HME}.
\end{Remark}

\begin{Remark}
Assume the set-up of Proposition~\ref{prop:MEslice}. The assumption that $\cG_N$ is source connected is quite strong and has the following consequences.
\begin{enumerate}\itemsep=0pt
\item[(i)] Denote by $p$ the unique point such that $N\cap L=\{p\}$.
Then the isotropy group $\cG_p^p$ of $\cG$ at~$p$ is connected. Indeed, it equals \smash{$(\cG_N)_p^p=\bt_{\cG_N}^{-1}(p)$} and $\cG_N$ is source (and target) connected by assumption. It follows that the isotropy group of the holonomy groupoid of $\cF$ at $p$ is connected, since the holonomy groupoid
 is a quotient of $\cG$ \cite[Example~3.4\,(4)]{AndrSk}. This is consistent with Remark~\ref{rem:Hppconn}.

\item[(ii)] The intersections of the leaves of $\cF$ with $N$ are connected. Indeed, they agree with the leaves of $\cG_N$, and the latter is source connected by assumption.
\end{enumerate}
\end{Remark}

\begin{Example}
Let \begin{itemize}\itemsep=0pt
\item
$V\to M$ be a vector bundle,
\item $G$ be a connected Lie group acting equivariantly on $V$ (not necessarily by vector bundle automorphisms) such that the action preserves the zero section $M$ and is transitive there.
\end{itemize}
Denote by $\cF$ the singular foliation on $V$ given by the infinitesimal generators of the action (hence~$M$ is a leaf).
Let $p\in M$, so the fiber $N:=V_p$ is a slice of the singular foliation through~$p$. If the isotropy group $G_p$ at $p$ is connected, then
$\bigl(N,\iota_N^{-1}\cF\bigr)$ and $(V,\cF)$ are Hausdorff Morita equivalent singular foliations.

Indeed, a source connected Lie groupoid giving rise to $\cF$ is the transformation groupoid $\cG=G\ltimes V$. The restricted Lie groupoid $\cG_N$ is $G_p\ltimes N$, the transformation groupoid of the $G_p$ action on $N$; by assumption it is source connected, so we can apply Proposition~\ref{prop:MEslice}.

A special case (where all the foliations involved are linear) is obtained starting with a Lie group $G$ and a connected closed Lie subgroup $H$. The action of $G$ on $M:=G/H$ lifts to an~equivariant linear action of $G$ on the tangent bundle of $M$ (and on any associated bundle). This example was inspired by \cite[Example~3.9]{CamilleSR}, which corresponds to the case $G={\rm SO}(3)$ and $H={\rm SO}(2)$ (so $\cF$ is a singular foliation on the tangent bundle of $M=S^2$, arising from the~${\rm SO}(3)$ action by rotations).
\end{Example}

Theorem~\ref{thm:main} and Proposition~\ref{prop:MEslice} imply the following.

\begin{Corollary}\label{cor:slicelin}
Let $\cF$ be a singular foliation arising from a Hausdorff, source connected, Lie groupoid $\cG$.
Let $L$ be an embedded leaf, $N$ a slice. Suppose the restricted Lie groupoid $\cG_N:=\bt^{-1}(N)\cap \bs^{-1}(N)$ is source connected.
Then, whenever $\bigl(N,\iota_N^{-1}\cF\bigr)$ is linearizable {around the point $N\cap L$}, the singular foliation $\cF$ is linearizable {around $L$}.
 \end{Corollary}

\subsection*{Acknowledgements}

M.Z.~thanks Camille Laurent-Gengoux for fruitful discussions and pointing out the reference~\cite{CerveauSing}, Nicolas Tribouillard for comments on this manuscript, and Alfonso Garmendia, Karandeep Singh, Llohann Speran\c{c}a, Joel Villatoro and Aldo Witte for useful discussions. Karandeep Singh also contributed directly to the proof of Proposition~\ref{prop:Cerveau}. We thank both referees for sharing insights that improved the paper; in particular, some of the arguments involving deformation spaces in Section~\ref{subsec:tne}, Remarks~\ref{rem:reffunct} and~\ref{rem:refhom} are due to them.
M.Z.~acknowledges partial support by IAP Dygest,
 Methusalem grants METH/15/026 and METH/21/03~-- long term structural funding of the Flemish Government, the FNRS and FWO under EOS projects G0H4518N and G0I2222N, the FWO research projects G083118N, G0B3523N and G014726N (Belgium).

\pdfbookmark[1]{References}{ref}
\LastPageEnding


\begin{thebibliography}{99}
\footnotesize\itemsep=0pt

\bibitem{AndrSk}
Androulidakis I., Skandalis G., The holonomy groupoid of a~singular foliation,
 \href{https://doi.org/10.1515/CRELLE.2009.001}{\textit{J.~Reine Angew.
 Math.}} \textbf{626} (2009), 1--37,
 \href{http://arxiv.org/abs/math.DG/0612370}{arXiv:math.DG/0612370}.

\bibitem{AZ1}
Androulidakis I., Zambon M., Smoothness of holonomy covers for singular
 foliations and essential isotropy,
 \href{https://doi.org/10.1007/s00209-013-1166-5}{\textit{Math.~Z.}}
 \textbf{275} (2013), 921--951,
 \href{http://arxiv.org/abs/1111.1327}{arXiv:1111.1327}.

\bibitem{AZ2}
Androulidakis I., Zambon M., Holonomy transformations for singular foliations,
 \href{https://doi.org/10.1016/j.aim.2014.02.003}{\textit{Adv. Math.}}
 \textbf{256} (2014), 348--397,
 \href{http://arxiv.org/abs/1205.6008}{arXiv:1205.6008}.

\bibitem{AZ6}
Androulidakis I., Zambon M., Stefan--{S}ussmann singular foliations, singular
 subalgebroids and their associated sheaves,
 \href{https://doi.org/10.1142/S0219887816410012}{\textit{Int.~J. Geom.
 Methods Mod. Phys.}} \textbf{13} (2016), 1641001, 17~pagges.

\bibitem{DefSpacesEulerLike}
Bischoff F., Bursztyn H., Lima H., Meinrenken E., Deformation spaces and normal
 forms around transversals,
 \href{https://doi.org/10.1112/s0010437x1900784x}{\textit{Compos. Math.}}
 \textbf{156} (2020), 697--732,
 \href{http://arxiv.org/abs/1807.11153}{arXiv:1807.11153}.

\bibitem{SplittingThmEulerlike}
Bursztyn H., Lima H., Meinrenken E., Splitting theorems for {P}oisson and
 related structures,
 \href{https://doi.org/10.1515/crelle-2017-0014}{\textit{J.~Reine Angew.
 Math.}} \textbf{754} (2019), 281--312,
 \href{http://arxiv.org/abs/1605.05386}{arXiv:1605.05386}.

\bibitem{CerveauSing}
Cerveau D., Distributions involutives singuli\`eres,
 \href{https://doi.org/10.5802/aif.761}{\textit{Ann. Inst. Fourier
 (Grenoble)}} \textbf{29} (1979), 261--294.

\bibitem{LinProperGroupoids}
Crainic M., Struchiner I., On the linearization theorem for proper {L}ie
 groupoids, \href{https://doi.org/10.24033/asens.2200}{\textit{Ann. Sci. \'Ec.
 Norm. Sup\'er.}} \textbf{46} (2013), 723--746,
 \href{http://arxiv.org/abs/1103.5245}{arXiv:1103.5245}.

\bibitem{MatiasME}
del Hoyo M.L., Lie groupoids and their orbispaces,
 \href{https://doi.org/10.4171/PM/1930}{\textit{Port. Math.}} \textbf{70}
 (2013), 161--209, \href{http://arxiv.org/abs/1212.6714}{arXiv:1212.6714}.

\bibitem{CamilleSR}
Fischer S.-R., Laurent-Gengoux C., A~classification of neighborhoods around
 leaves of a~singular foliation,
 \href{http://arxiv.org/abs/2401.05966}{arXiv:2401.05966}.

\bibitem{francis2024singularfoliationstangentgiven}
Francis M., On singular foliations tangent to a~given hypersurface,
 \href{https://doi.org/10.4171/JNCG/600}{\textit{J.~Noncommut. Geom.}}, {t}o
 appear, \href{http://arxiv.org/abs/2311.03940}{arXiv:2311.03940}.

\bibitem{AlfonsoThesis}
Garmendia A., Groupoids and singular foliations, Ph.D.~Thesis, {K}U Leuven, 2019,
 available at
 \url{https://perswww.kuleuven.be/~u0096206/students/thesisAlfonsoFInal.pdf},
 \href{http://arxiv.org/abs/2107.10502}{arXiv:2107.10502}.

\bibitem{OriAlfonso}
Garmendia A., Yudilevich O., On the inner automorphisms of a~singular
 foliation,
 \href{https://doi.org/10.1007/s00209-018-2212-0}{\textit{Math.~Z.}}
 \textbf{293} (2019), 725--729,
 \href{http://arxiv.org/abs/1804.06103}{arXiv:1804.06103}.

\bibitem{HME}
Garmendia A., Zambon M., Hausdorff {M}orita equivalence of singular foliations,
 \href{https://doi.org/10.1007/s10455-018-9620-6}{\textit{Ann. Global Anal.
 Geom.}} \textbf{55} (2019), 99--132,
 \href{http://arxiv.org/abs/1803.00896}{arXiv:1803.00896}.

\bibitem{HigsonReza}
Haj Saeedi~Sadegh A.R., Higson N., Euler-like vector fields, deformation spaces
 and manifolds with filtered structure,
 \href{https://doi.org/10.4171/DM/619}{\textit{Doc. Math.}} \textbf{23}
 (2018), 293--325, \href{http://arxiv.org/abs/1611.05312}{arXiv:1611.05312}.

\bibitem{singfolnotes}
Laurent-Gengoux C., Louis R., Ryvkin L., An invitation to singular foliations,
 {A}dvanced Courses in Mathematics -- CRM Barcelona, 2024,
 \href{http://arxiv.org/abs/2407.14932}{arXiv:2407.14932}.

\bibitem{NbhdSingLeaf}
Laurent-Gengoux C., Ryvkin L., The neighborhood of a~singular leaf,
 \href{https://doi.org/10.5802/jep.165}{\textit{J.~\'Ec. Polytech. Math.}}
 \textbf{8} (2021), 1037--1064,
 \href{http://arxiv.org/abs/2004.07019}{arXiv:2004.07019}.

\bibitem{symsingfol}
Louis R., On symmetries of singular foliations,
 \href{https://doi.org/10.1016/j.geomphys.2023.104833}{\textit{J.~Geom.
 Phys.}} \textbf{189} (2023), 104833, 31~pages,
 \href{http://arxiv.org/abs/2203.01585}{arXiv:2203.01585}.

\bibitem{MackenzieGrdAld}
Mackenzie K.C.H., General theory of {L}ie groupoids and {L}ie algebroids,
 \textit{London Math. Soc. Lecture Note Ser.}, Vol.~213,
 \href{https://doi.org/10.1017/CBO9781107325883}{Cambridge University Press},
 Cambridge, 2005.

\bibitem{ELisotropic}
Meinrenken E., Euler-like vector fields, normal forms, and isotropic
 embeddings, \href{https://doi.org/10.1016/j.indag.2020.08.006}{\textit{Indag.
 Math.~(N.S.)}} \textbf{32} (2021), 224--245,
 \href{http://arxiv.org/abs/2001.10518}{arXiv:2001.10518}.

\bibitem{MarcutThesis}
M\u{a}rcu\c{t} I., Normal forms in Poisson geometry, Ph.D.~Thesis, {U}trecht
 University, 2013, \href{http://arxiv.org/abs/1301.4571}{arXiv:1301.4571}.

\bibitem{WitteCohomElliptic}
Witte A., The cohomology of the elliptic tangent bundle,
 \href{https://doi.org/10.1016/j.indag.2021.09.003}{\textit{Indag. Math.
 (N.S.)}} \textbf{33} (2022), 372--387,
 \href{http://arxiv.org/abs/2104.04845}{arXiv:2104.04845}.

\end{thebibliography}
\end{document}